\documentclass[12pt]{article}

\usepackage{amsmath, epsfig, cite}
\usepackage{amssymb}
\usepackage{amsfonts}
\usepackage{latexsym}

\newtheorem{thm}{Theorem}[section]

\newtheorem{conj}[thm]{Conjecture}
\newtheorem{lem}[thm]{Lemma}
\newtheorem{prob}[thm]{\it Problem}
\newcommand{\pf}{\noindent{\it Proof.} }
\newcommand*{\fl}[2]{\left\lfloor\frac{#1}{#2}\right\rfloor}
\newcommand*{\pp}[2]{ \left\{\frac{#1}{#2} \right\}}

\numberwithin{equation}{section}

\newcommand{\qed}{{\hfill$\square$}\medskip}

\begin{document}

\nocite{*}
\begin{center}
{\Large\bf Proof of a conjecture of Z.-W. Sun on\\[5pt] the divisibility of a triple sum}
\end{center}

\vskip 2mm \centerline{Victor J. W. Guo$^1$ and Ji-Cai Liu$^2$}
\begin{center}
{\footnotesize Department of Mathematics, Shanghai Key Laboratory of
PMMP, East China Normal University,\\ 500 Dongchuan Road, Shanghai
200241,
 People's Republic of China\\
$^1${\tt jwguo@math.ecnu.edu.cn,\quad
http://math.ecnu.edu.cn/\textasciitilde{jwguo}}  \quad {$^2$\tt jc2051@163.com} }
\end{center}


\vskip 0.7cm \noindent{\bf Abstract.} The numbers $R_n$ and $W_n$ are defined as
\begin{align*}
R_n=\sum_{k=0}^{n}{n+k\choose 2k}{2k\choose k}\frac{1}{2k-1},\ \text{and}\ W_n=\sum_{k=0}^{n}{n+k\choose 2k}{2k\choose k}\frac{3}{2k-3}.
\end{align*}
We prove that, for any positive integer $n$ and odd prime $p$, there hold
\begin{align*}
\sum_{k=0}^{n-1}(2k+1)R_k^2 &\equiv 0 \pmod{n}, \\
\sum_{k=0}^{p-1}(2k+1)R_k^2 &\equiv 4p(-1)^{\frac{p-1}{2}} -p^2 \pmod{p^3}, \\
9\sum_{k=0}^{n-1}(2k+1)W_k^2 &\equiv 0 \pmod{n},  \\
\sum_{k=0}^{p-1}(2k+1)W_k^2  &\equiv  12p(-1)^{\frac{p-1}{2}}-17p^2 \pmod{p^3}, \quad\text{if $p>3$.}
\end{align*}
The first two congruences were originally conjectured by Z.-W. Sun. Our proof is based on the multi-variable Zeilberger algorithm
and the following observation:
$$
{2n\choose n}{n\choose k}{m\choose k}{k\choose m-n}\equiv 0\pmod{{2k\choose k}{2m-2k\choose m-k}},
$$
where $0\leqslant k\leqslant n\leqslant m \leqslant 2n$.


\vskip 2mm
\noindent{\it MR Subject Classifications}: 11A07, 11B65, 05A10

\section{Introduction}
For any nonnegative integer $n$, it is easy to see that ${2n\choose n}\frac{1}{2n-1}$ is always an integer.
Z.-W. Sun \cite{Sun} introduced the following numbers
\begin{align*}
R_n=\sum_{k=0}^{n}{n+k\choose 2k}{2k\choose k}\frac{1}{2k-1}
\end{align*}
and proved some interesting arithmetic properties of these numbers. For example, Sun \cite{Sun} proved that, if $p$ is a prime of the form $4k+1$, then
$$
R_{\frac{p-1}{2}}\equiv p-(-1)^{\frac{p-1}{4}} 2x \pmod{p^2},
$$
where $p=x^2+y^2$ with $x\equiv 1\pmod{4}$.

The first aim of this paper is to prove the following result, which was originally conjectured by Z.-W. Sun (see \cite[Conjecture 5.4]{Sun}).
\begin{thm}\label{thm:2}
Let $n$ be a positive integer and $p$ an odd prime. Then
\begin{align}
\sum_{k=0}^{n-1}(2k+1)R_k^2 &\equiv 0 \pmod{n}, \label{eq:5} \\
\sum_{k=0}^{p-1}(2k+1)R_k^2 &\equiv 4p(-1)^{\frac{p-1}{2}} -p^2 \pmod{p^3}. \label{eq:new5}
\end{align}
\end{thm}

Since
\begin{align}
{2n\choose n}\frac{3}{2n-3}={2n\choose n}\frac{1}{2n-1}+{2n-2\choose n-1}\frac{8}{2n-3}, \label{eq:2k-3}
\end{align}
we see that ${2n\choose n}\frac{3}{2n-3}$ is always an integer.

Let
\begin{align*}
W_n=\sum_{k=0}^{n}{n+k\choose 2k}{2k\choose k}\frac{3}{2k-3}.
\end{align*}
The second aim of this paper is to prove the following congruence and supercongruence.
\begin{thm}\label{thm:3}
Let $n$ be a positive integer and let $p>3$ be a prime. Then
\begin{align}
9\sum_{k=0}^{n-1}(2k+1)W_k^2 &\equiv 0 \pmod{n}, \label{eq:6} \\
\sum_{k=0}^{p-1}(2k+1)W_k^2  &\equiv  12p(-1)^{\frac{p-1}{2}}-17p^2 \pmod{p^3}.    \label{eq:new6}
\end{align}
\end{thm}

The rest of the paper is organized as follows. In the next section, we shall establish some important lemmas,
including a result on the divisibility of a product of four binomial coefficients.
The proofs of Theorems \ref{thm:2} and \ref{thm:3} will be given in Sections 3 and 4, respectively.
In Section 5, we will propose some open problems for further study.

\section{Some lemmas}

We first establish the following lemma.
\begin{lem}\label{lem:one}
Let $k,m,n$ be nonnegative integers with $0\leqslant k\leqslant n\leqslant m \leqslant 2n$. Then
\begin{align}
{2n\choose n}{n\choose k}{m\choose k}{k\choose m-n}\equiv 0\pmod{{2k\choose k}{2m-2k\choose m-k}}.\label{eq:4}
\end{align}
\end{lem}

\noindent{\it Proof.} Let $p$ be a prime. For the $p$-adic order of $n!$, there is a famous formula
$$\text{ord}_p n!=\sum_{j=1}^{\infty}\fl{n}{p^j},$$
where $\left\lfloor x\right\rfloor$ denotes the greatest integer less than or equal to $x$.
It follows that
\begin{align}
&\hskip -2mm\text{ord}_p \frac{{2n\choose n}{n\choose k}{m\choose k}{k\choose m-n}}{{2k\choose k}{2m-2k\choose m-k}} \nonumber \\
&=\text{ord}_p \frac{(2n)!(m)!(m-k)!k!}{(n-k)!(k-m+n)!(2m-2k)!(2k)!n!(m-n)!}  \nonumber \\
&=\sum_{j=0}^{\infty}\Bigg( \fl{2n}{p^j}+\fl{m}{p^j}+\fl{m-k}{p^j}+\fl{k}{p^j} \nonumber \\
&\quad{}-\fl{n-k}{p^j}-\fl{k-m+n}{p^j}-\fl{2m-2k}{p^j}-\fl{2k}{p^j}-\fl{n}{p^j}-\fl{m-n}{p^j}\Bigg).  \label{eq:factorfour}
\end{align}
To prove \eqref{eq:4}, it suffices to show that the right-hand side of \eqref{eq:factorfour} is nonnegative for any prime $p$.
In what follows, we shall prove that, for any positive integer $d$, there holds
\begin{align}
&\fl{2n}{d}+\fl{m}{d}+\fl{m-k}{d}+\fl{k}{d} \nonumber\\
&\geqslant \fl{n-k}{d}+\fl{k-m+n}{d}+\fl{2m-2k}{d}+\fl{2k}{d}+\fl{n}{d}+\fl{m-n}{d}.  \label{eq:fourcases}
\end{align}
Note that, by \cite[Theorem 1.3]{Bober}, we have
\begin{align}
\fl{2n}{d}+\fl{k}{d} \geqslant \fl{n}{d}+\fl{n-k}{d}+\fl{2k}{d}. \label{eq:1}
\end{align}

We proceed by considering three cases.
\begin{itemize}
\item If $\fl{2m-2k}{d}=2\fl{m-k}{d}$, then by \eqref{eq:1} and the fact $\fl{m}{d}\geqslant \fl{k-m+n}{d}+\fl{m-k}{d}+\fl{m-n}{d}$
the inequality \eqref{eq:fourcases} holds.

\item If $\fl{2m-2k}{d}=2\fl{m-k}{d}+1$ and $\fl{m}{d}\geqslant \fl{k-m+n}{d}+\fl{m-k}{d}+\fl{m-n}{d}+1$, then by \eqref{eq:1} again,
the inequality \eqref{eq:fourcases} holds.

\item If $\fl{2m-2k}{d}=2\fl{m-k}{d}+1$ and $\fl{m}{d}=\fl{k-m+n}{d}+\fl{m-k}{d}+\fl{m-n}{d}$, then $\pp{m-k}{d}\geqslant\frac{1}{2}$ and
$\pp{k-m+n}{d}+\pp{m-k}{d}+\pp{m-n}{d}<1$, where $\{x\}=x-\lfloor x\rfloor$ denotes the fraction part of $x$.
It follows that
\begin{align*}
\pp{n}{d}&=\pp{k-m+n}{d}+\pp{m-k}{d}\geqslant \frac{1}{2}, \\
\pp{k}{d}&=\pp{k-m+n}{d}+\pp{m-n}{d}<\frac{1}{2}.
\end{align*}
Hence, we have $\fl{2n}{d}=2\fl{n}{d}+1$ and $\fl{2k}{d}=2\fl{k}{d}$. Since $\fl{n}{d}\geqslant \fl{n-k}{d}+\fl{k}{d}$,
we see that the inequality \eqref{eq:fourcases} still holds in this case.
\end{itemize}
This completes the proof.    \qed

\noindent{\it Remark.} Some similar divisibility properties of products of binomial coefficients have already been
obtained by Z.-W. Sun \cite{Sun2012,Sun2013}, Guo \cite{Guo2013,Guo2014}, and Guo and Krattenthaler \cite{GK}.

\medskip
The following lemma is critical in our proof of \eqref{eq:5}.
\begin{lem}\label{lem:3} Let $n$ be a nonnegative integer. Then
\begin{align}
{2n\choose n}\sum_{m=n}^{2n}\sum_{k=0}^n{n\choose k}{m\choose k}{k\choose m-n}\frac{1}{(2k-1)(2m-2k-1)}\equiv 0\pmod{n+1}. \label{eq:10}
\end{align}
\end{lem}
\pf By Lemma \ref{lem:one} and the fact $(2k-1)\mid{2k\choose k}$, we deduce that the left-hand side of \eqref{eq:10},
denoted by $X_n$, is always an integer.
Applying the multi-variable Zeilberger algorithm (see \cite{AZ,WZ}), we find that the numbers $X_n$ satisfy the following fifth-order recurrence:
\begin{align}
&512(n-1)(n+2)(2n-1)(2n+1)(2n+3)^2(7056n^5+90958n^4+438943n^3+960044n^2 \nonumber\\
&+877175n+187500)X_n-64(n+1)(2n+1)(2n+3)(7056n^8+6286n^7-938209n^6\nonumber\\
&-8101907n^5-29947351n^4-56721017n^3-53422948n^2-18893910n+1125000)X_{n+1}\nonumber\\
&-16(n+1)(n+2)(2n+3)(225792n^8+3729152n^7+25149436n^6+88143056n^5 \nonumber\\
&+165692905n^4+141655826n^3-9994217n^2-100176150n-47709000)X_{n+2} \nonumber\\
&+4(n+1)(n+2)(n+3)(56448n^8+558320n^7+373360n^6-15659602n^5-82475561n^4\nonumber\\
&-185781850n^3-205165415n^2-98261100n-10317600)X_{n+3}\nonumber\\
&+4(n+1)(n+2)(n+3)(n+4)^2(56448n^6+628880n^5+2620038n^4+4826445n^3\nonumber\\
&+2994664n^2-1570545n-1935450)X_{n+4}-(n+1)(n+2)(n+3)(n+4)(n+5)^2\nonumber \\
&\times(7056n^5 +55678n^4+145671n^3+118403n^2-54636n-84672)X_{n+5}=0. \label{eq:11}
\end{align}
It is interesting that we can deduce the following third-order recurrence for $X_n$ from \eqref{eq:11}:
\begin{align}
&\hskip -2mm
128(n-1)(n+2)(2n+1)(2n+3)(7n+15)X_{n} \nonumber\\
&-16(n+1)(7n^4-6n^3-121n^2-210n-90)X_{n+1} \nonumber\\
&-4(n+1)( n+2)^2( 56n^2+127n+57 )X_{n+2} \nonumber \\
&+(n+1)(n+2)(n+3)^2(7n+8)X_{n+3}=0. \label{eq:12}
\end{align}
In fact, if we denote the left-hand sides of \eqref{eq:11} and \eqref{eq:12} by $\alpha_n$ and $\beta_n$, respectively,
then we can easily check that
\begin{align*}
&(7n+15 )(7n+22)\alpha_n-4(2n-1)(2n+3)(7n+22)(7056n^5+90958n^4+438943n^3 \nonumber\\
&+960044n^2+877175n+187500)\beta_n -16(n+1)(12348n^6+175028n^5+969283n^4\nonumber\\
&+2677654n^3+3865514n^2+2712498n+679905)\beta_{n+1}+(n+1)(n+2)(7n+15 )\nonumber\\
&\times(7056n^5+55678n^4+145671n^3+118403n^2-54636n-84672)\beta_{n+2}=0.
\end{align*}
Therefore, by induction on $n$, we immediately obtain $\beta_n=0$, i.e., the recurrence \eqref{eq:12} is true.
It follows from \eqref{eq:12} that $(n+1)$ divides
$$
128(n-1)(n+2)(2n+1)(2n+3)(7n+15)X_{n}.
$$
Since $\gcd(n+1,(n+2)(2n+1)(2n+3))=1$, $\gcd(n+1,n-1)=1,2$, and $\gcd(n+1,7n+15)=\gcd(n+1,8)=1,2,4,8$.
we see that $(n+1)$ divides $2^{11} X_{n}$. Namely, the expression $\frac{2^{11}X_n}{n+1}$ is an integer,
which means that $\frac{X_n}{n+1}$ must be integral,
since the denominator of $\frac{X_n}{n+1}$ is odd. This completes the proof. \qed

The following lemma plays an important part in our proof of \eqref{eq:6}.
\begin{lem}\label{lem:4}
Let $n$ be a nonnegative integer. Then
\begin{align}
{2n\choose n}\sum_{m=n}^{2n}\sum_{k=0}^n{n\choose k}{m\choose k}{k\choose m-n}\frac{81}{(2k-3)(2m-2k-3)}\equiv 0\pmod{n+1}. \label{eq:16}
\end{align}
\end{lem}
\pf Let $9Y_n$ denote the left-hand side of \eqref{eq:16}. By Lemma \ref{lem:one} and \eqref{eq:2k-3}, we know that
$Y_n$ is an integer. Applying the multi-variable Zeilberger algorithm, we obtain
\begin{align}
&-512(n-3)(n+2)(2n-5)(2n+1)(2n+3)^2(2n+5)(2n+7)(2n+9)(3993696n^{10} \nonumber\\
&+101741444n^9+1106902594n^8+6675903296n^7+24039265882n^6+50526147407n^5 \nonumber\\
&+49093431499n^4-22567478757n^3-115591006351n^2-118410894910n   \nonumber\\
&-43001171400)Y_n+64(n+1)(2n+3)(2n+5)(2n+7)(2n+9)(7987392n^{14} \nonumber\\
&-175918232n^{13}-7078653400n^{12}-80904072538n^{11}-450205035754n^{10} \nonumber\\
&-1235909690096n^9-684324942173n^8+5752795629096n^7+17070643824448n^6  \nonumber\\
&+16426830143582n^5-8252471533811n^4-32840781231384n^3-29184924754630n^2  \nonumber\\
&-11229082751700n-1771825878000) Y_{n+1}+16(n+1)(n+2)(2n+5)(2n+7) \nonumber\\
&\times (2n+9)( 255596544n^{14}+6351704576n^{13}+65334850448n^{12}+346171943648n^{11}\nonumber\\
&+870630300008n^{10}-154616465854n^9-7786336958932n^8-23352735092682n^7 \nonumber\\
&-30879971643605n^6-8340607431055n^5+32770026974177n^4+51322013010691n^3 \nonumber\\
&+35260464547296n^2+12208208518740n+1725138622800) Y_{n+2}-4(n+1)(n+2) \nonumber\\
&\times (n+3)(2n+3)(2n+7)(2n+9)(31949568n^{13}+47141920n^{12}-8899203264n^{11} \nonumber\\
&-119872488660n^{10}-729008210810n^9-2427338831964n^8-4192616628250n^7   \nonumber\\
&-1193459139415n^6+10701703824509n^5+24849269008557n^4+27637389771751n^3     \nonumber\\
&+17335661249538n^2+5872316424120n+832457390400)Y_{n+3}    \nonumber\\
&-4(n+1)(n+2)(n+3)(n+4)(2n+3)(2n+5)(2n+9)(31949568n^{12}+710095456n^{11}  \nonumber\\
&+6681317524n^{10}+34029339728n^9+96227598957n^8+112933571943n^7    \nonumber\\
&-172413925739n^6-927631519653n^5-1726986013160n^4-1811282710094n^3   \nonumber\\
&-1125028801230n^2-387378800500n-56652486000) Y_{n+4}+(n+1)(n+2)   \nonumber\\
&\times (n+3)(n+4)(n+5)^2(2n+3)(2n+5)(2n+7)(3993696n^{10}+61804484n^9  \nonumber \\
&+370945918n^8+1004131008n^7+593610306n^6-3689013381n^5-11102496870n^4   \nonumber\\
&-14727940451n^3-10684139174n^2-4158064376n-674002560)Y_{n+5}=0. \label{eq:18}
\end{align}
Similarly as before, we can deduce the following simpler recurrence for $Y_n$ from \eqref{eq:18}:
\begin{align}
&-128(n-3)(n+2)(2n+1)(2n+3)^2(2n+5)(63n^3+390n^2+785n+506)Y_{n}    \nonumber \\
& +16(n+1)(2n+3)(2n+5) (63n^6-933n^5-7645n^4-17421n^3-13730n^2-2538n     \nonumber\\
&-252)Y_{n+1} +4(n+1)(n+2)(2n+1)(2n+5)(504n^5+2805n^4+5464n^3+4575n^2\nonumber\\
&+1400n-60)Y_{n+2} -(n+1)(n+2)(n+3)^2(2n+1) (2n+3) (63n^3+201n^2+194n\nonumber\\
&+48) Y_{n+3}=0, \label{eq:19}
\end{align}
by noticing that
\begin{align*}
&(63n^3+390n^2+785n+506)(63n^3+579n^2+1754n+1744)\gamma_n-4(2n-5)(2n+7) \nonumber\\
&\times (2n+9)(63n^3+579n^2+1754n+1744)( 3993696n^{10}+101741444n^9+1106902594n^8   \nonumber\\
&+6675903296n^7+24039265882n^6+50526147407n^5+49093431499n^4-22567478757n^3\nonumber\\
&-115591006351n^2-118410894910n-43001171400)\delta_n-4( n+1)(2n+9 )(754808544n^{13}\nonumber\\
&+21322656936n^{12}+262462697910n^{11}+1845558281063n^{10}+8129372080496n^9\nonumber\\
&+22891046730211n^8+38943827465846n^7+28214897181357n^6-31234241796612n^5     \nonumber\\
&-104745740003975n^4-123454281828448n^3-77917639095288n^2-25845991472440n \nonumber\\
&-3509843409600)\delta_{n+1} +( n+1)( n+2) (2n+3)(63n^3+390n^2+785n+506 )(3993696n^{10}    \nonumber\\
&+61804484n^9+370945918n^8+1004131008n^7+593610306n^6-3689013381n^5\nonumber\\
&-11102496870n^4-14727940451n^3-10684139174n^2-4158064376n-674002560)\delta_{n+3} \nonumber \\
&=0,
\end{align*}
where $\gamma_n$ and $\delta_n$ denote the left-hand sides of
\eqref{eq:18} and \eqref{eq:19}, respectively. From \eqref{eq:19},
it is easy to see that $(n+1)$ divides
$$
128(n-3)(n+2)(2n+1)(2n+3)^2(2n+5)(63n^3+390n^2+785n+506)Y_{n}.
$$
Since $\gcd(n+1,(n+2)(2n+1)(2n+3)^2)=1$, $\gcd(n+1,n-3)=2^i$ with $0\leqslant i\leqslant 2$,
$\gcd(2n+5,n+1)=3^j$ with $0\leqslant j\leqslant 1$, and $\gcd(n+1,63n^3+390n^2+785n+506)=\gcd(n+1,48)=2^u\cdot 3^v$
with $0\leqslant u\leqslant 4$ and $0\leqslant v\leqslant 1$, we conclude that
$(n+1)$ divides $2^{13}\cdot 9Y_{n}$. In other words, the ratio $\frac{2^{13}\cdot 9Y_n}{n+1}$ is an integer,
which means that $\frac{9Y_n}{n+1}$ must be integral, for the denominator of $\frac{9Y_n}{n+1}$ is odd.  This completes the proof.
\qed

To prove the supercongruence \eqref{eq:new5}, we further need the following lemma.
\begin{lem}Let $p$ be a prime. Then
\begin{align}
\sum_{i=0}^{p-1}\sum_{j=0}^{p-1}\frac{p^2(-1)^{i+j}}{(2i-1)(2j-1)(i+j+1)}
&\equiv 4p(-1)^{\frac{p-1}{2}} +3p^2 \pmod{p^3} ,  \label{eq:p2i-1}
\end{align}
\end{lem}
\noindent{\it Proof.}
It is easy to see that
\begin{align}
&\hskip -2mm
\sum_{i=0}^{p-1}\sum_{j=0}^{p-1}\frac{p^2(-1)^{i+j}}{(2i-1)(2j-1)(i+j+1)}  \nonumber \\
&=\sum_{i=0}^{p-1}\sum_{j=0}^{p-1}\left(\frac{p^2(-1)^{i+j}}{2(2i-1)(2j-1)}
-\frac{p^2(-1)^{i+j}}{4(2i-1)(i+j+1)}-\frac{p^2(-1)^{i+j}}{4(2j-1)(i+j+1)}\right) \nonumber \\
&=\frac{1}{2}\left(\sum_{i=0}^{p-1}\frac{p(-1)^{i}}{2i-1}\right)^2
-\frac{1}{2}\sum_{i=0}^{p-1}\sum_{j=0}^{p-1}\frac{p^2(-1)^{i+j}}{(2i-1)(i+j+1)}. \label{eq:sum-part1}
\end{align}
Since
\begin{align}
\sum_{i=0}^{p-1}\frac{p(-1)^{i}}{2i-1}
&=-2p+(-1)^{\frac{p+1}{2}}+\sum_{k=1}^{\frac{p-3}{2}}(-1)^{\frac{p+1}{2}+k}\left(\frac{p}{p-2k}+\frac{p}{p+2k}\right) \notag  \\
&=-2p+(-1)^{\frac{p+1}{2}}+\sum_{k=1}^{\frac{p-3}{2}}(-1)^{\frac{p+1}{2}+k}\frac{2p^2}{p^2-4k^2} \notag \\
&\equiv -2p+(-1)^{\frac{p+1}{2}}+\frac{1}{2}\sum_{k=1}^{\frac{p-3}{2}}(-1)^{\frac{p-1}{2}+k}\frac{p^2}{k^2}
\pmod{p^3},  \label{eq:sum-part200}
\end{align}
we have
\begin{align}
\left(\sum_{i=0}^{p-1}\frac{p(-1)^{i}}{2i-1}\right)^2
\equiv 1+4p(-1)^{\frac{p-1}{2}}+4p^2-p^2\sum_{k=1}^{\frac{p-3}{2}}\frac{(-1)^k}{k^2} \pmod{p^3}.  \label{eq:sum-part2}
\end{align}
In view of \eqref{eq:sum-part1} and \eqref{eq:sum-part2}, to prove \eqref{eq:p2i-1}, it suffices to prove that
\begin{align}
\sum_{i=0}^{p-1}\sum_{j=0}^{p-1}\frac{p^2(-1)^{i+j}}{(2i-1)(i+j+1)}
\equiv 1-4p(-1)^{\frac{p-1}{2}}-2p^2-p^2\sum_{k=1}^{\frac{p-3}{2}}\frac{(-1)^k}{k^2} \label{eq:sum-part3}
\pmod{p^3}.
\end{align}

It is clear that
\begin{align}
\sum_{i=0}^{p-1}\sum_{j=0}^{p-1}\frac{p^2(-1)^{i+j}}{(2i-1)(i+j+1)}=\sum_{i=0}^{p-1}\sum_{j=0}^{p-1}\left(\frac{2p^2(-1)^{i+j}}{(2i-1)(2j+3)}
-\frac{p^2(-1)^{i+j}}{(i+j+1)(2j+3)}\right), \label{eq:identity}
\end{align}
and
\begin{align}
&\sum_{i=0}^{p-1}\sum_{j=0}^{p-1}\frac{p^2(-1)^{i+j}}{(2i-1)(2j+3)} \notag\\
&=\sum_{i=0}^{p-1}\sum_{j=2}^{p+1}\frac{p^2(-1)^{i+j}}{(2i-1)(2j-1)} \notag\\
&=\sum_{i=0}^{p-1}\sum_{j=0}^{p-1}\frac{p^2(-1)^{i+j}}{(2i-1)(2j-1)}+\sum_{i=0}^{p-1}\left(\frac{2p^2(-1)^{i}}{2i-1}
-\frac{p^2(-1)^{i}}{(2i-1)(2p-1)}+\frac{p^2(-1)^{i}}{(2i-1)(2p+1)}\right) \notag \\
&=\left(\sum_{i=0}^{p-1}\frac{p(-1)^{i}}{2i-1}\right)^2+\sum_{i=0}^{p-1}\left(\frac{2p^2(-1)^{i}}{2i-1}
-\frac{2p^2(-1)^{i}}{(4p^2-1)(2i-1)}\right)  \notag\\
&\equiv\left(\sum_{i=0}^{p-1}\frac{p(-1)^{i}}{2i-1}\right)^2+4\sum_{i=0}^{p-1}\frac{p^2(-1)^{i}}{2i-1} \notag\\
&\equiv 1-4p^2-p^2\sum_{k=1}^{\frac{p-3}{2}}\frac{(-1)^k}{k^2} \pmod{p^3} \quad\text{(by \eqref{eq:sum-part200} and \eqref{eq:sum-part2}).}
\label{eq:congr-two}
\end{align}
On the other hand, we have
\begin{align}
&\hskip -1mm \sum_{i=0}^{p-1}\sum_{j=0}^{p-1}\frac{p^2(-1)^{i+j}}{(i+j+1)(2j+3)} \notag \\
&=\sum_{i=-2}^{p-3}\sum_{j=2}^{p+1}\frac{p^2(-1)^{i+j}}{(i+j+1)(2j-1)} \notag \\
&=\sum_{i=0}^{p-1}\sum_{j=0}^{p-1}\frac{p^2(-1)^{i+j}}{(i+j+1)(2j-1)}  +\sum_{j=0}^{p-1}\left(\frac{1}{j+1}-\frac{1}{j+2}
+\frac{1}{j+p+1}-\frac{1}{j+p+2}\right)\frac{p^2(-1)^{j}}{2j+3} \notag\\
&\quad{}+\sum_{i=0}^{p-1}\left(\frac{p^2(-1)^{i}}{i+1}+\frac{p^2(-1)^{i}}{i+2}
-\frac{p^2(-1)^{i}}{(i+p+1)(2p-1)}+\frac{p^2(-1)^{i}}{(i+p+2)(2p+1)}\right).  \label{eq:sum-hard}
\end{align}
It is not hard to see that
\begin{align}
&\hskip -2mm\sum_{i=0}^{p-1}\left(\frac{p^2(-1)^{i}}{i+1}+\frac{p^2(-1)^{i}}{i+2}\right)=p^2+\frac{p^2}{p+1}\equiv 2p^2 \pmod{p^3},
\label{eq:sum-lem-1} \\[5pt]
&\hskip -2mm\sum_{i=0}^{p-1}\left(-\frac{p^2(-1)^{i}}{(i+p+1)(2p-1)}+\frac{p^2(-1)^{i}}{(i+p+2)(2p+1)}\right) \notag\\
&\equiv -\frac{p}{4p-2}-\frac{p}{4p+2}+p^2 +\sum_{i=0}^{p-2}\frac{p^2(-1)^{i}}{i+1}+ \sum_{i=0}^{p-3}\frac{p^2(-1)^{i}}{i+2}
 \notag \\
&=-\frac{2p^2}{4p^2-1}+2p^2 \notag \\
&\equiv 4p^2 \pmod{p^3},  \label{eq:sum-lem-2} \\[5pt]
&\hskip -2mm \sum_{i=0}^{p-1}\frac{p^2(-1)^{i}}{2i+3} =2p^2+\frac{p^2}{2p+1}-\frac{p^2}{2p-1}+\sum_{i=0}^{p-1}\frac{p^2(-1)^{i}}{2i-1}
\equiv 2p^2-p(-1)^{\frac{p-1}{2}} \pmod{p^3}, \notag
\end{align}
where we used \eqref{eq:sum-part200} in the last step, and so
\begin{align}
&\hskip -2mm \sum_{j=0}^{p-1}\left(\frac{1}{j+1}-\frac{1}{j+2}
+\frac{1}{j+p+1}-\frac{1}{j+p+2}\right)\frac{p^2(-1)^{j}}{2j+3}  \notag\\
&=\sum_{j=0}^{p-1}\left(\frac{p^2(-1)^{j}}{j+1}-\frac{4p^2(-1)^{j}}{2j+3}+\frac{p^2(-1)^{j}}{j+2}\right)
+\frac{1}{2p-1}\sum_{j=0}^{p-1}\left(\frac{2p^2(-1)^{j}}{2j+3}-\frac{p^2(-1)^{j}}{j+p+1}\right) \notag \\
&\quad{}+\frac{1}{2p+1}\sum_{j=0}^{p-1}\left(\frac{p^2(-1)^{j}}{j+p+2}-\frac{2p^2(-1)^{j}}{2j+3}\right) \notag \\
&=\sum_{j=0}^{p-1}\left(\frac{p^2(-1)^{j}}{j+1}+\frac{p^2(-1)^{j}}{j+2}\right)
+\left(-4+\frac{2}{2p-1}-\frac{2}{2p+1}\right)\sum_{j=0}^{p-1}\frac{p^2(-1)^{j}}{2j+3}  \notag\\
&\quad{}+\sum_{j=0}^{p-1}\left(\frac{p^2(-1)^{j}}{(2p+1)(j+p+2)}-\frac{p^2(-1)^{j}}{(2p-1)(j+p+1)}\right) \notag \\
&\equiv 2p^2+\left(-4+\frac{4}{4p^2-1}\right)(2p^2-p(-1)^{\frac{p-1}{2}})+4p^2  \notag \\
&\equiv -10p^2+8p(-1)^{\frac{p-1}{2}}\pmod{p^3}.   \label{eq:sum-lem-3}
\end{align}

Substituting \eqref{eq:sum-lem-1}--\eqref{eq:sum-lem-3} into \eqref{eq:sum-hard}, we have
\begin{align}
\sum_{i=0}^{p-1}\sum_{j=0}^{p-1}\frac{p^2(-1)^{i+j}}{(i+j+1)(2j+3)}
\equiv\sum_{i=0}^{p-1}\sum_{j=0}^{p-1}\frac{p^2(-1)^{i+j}}{(i+j+1)(2j-1)}
-4p^2+8p(-1)^{\frac{p-1}{2}}\pmod{p^3},  \label{eq:congr-last-1}
\end{align}
while substituting \eqref{eq:congr-two} into \eqref{eq:identity}, we get
\begin{align}
&\hskip -2mm \sum_{i=0}^{p-1}\sum_{j=0}^{p-1}\frac{p^2(-1)^{i+j}}{(i+j+1)(2j+3)}
+\sum_{i=0}^{p-1}\sum_{j=0}^{p-1}\frac{p^2(-1)^{i+j}}{(i+j+1)(2j-1)}  \notag \\
&\equiv 2-8p^2-2p^2\sum_{k=1}^{\frac{p-3}{2}}\frac{(-1)^k}{k^2} \pmod{p^3}. \label{eq:congr-last-2}
\end{align}
Combining \eqref{eq:congr-last-1} and \eqref{eq:congr-last-2}, and noticing the symmetry of $i$ and $j$, we immediately
obtain \eqref{eq:sum-part3}. This completes the proof. \qed

We now give our last lemma, which is related to \eqref{eq:new6}.

\begin{lem}Let $p>3$ be a prime. Then
\begin{align}
\sum_{i=0}^{p-1}\sum_{j=0}^{p-1}\frac{p^2(-1)^{i+j}}{(2i-3)(2j-3)(i+j+1)}
&\equiv \frac{4p}{3}(-1)^{\frac{p-1}{2}}+\frac{p^2}{3}
\pmod{p^3}.  \label{eq:p2i-3}
\end{align}
\end{lem}
\noindent{\it Proof.} Similarly as \eqref{eq:sum-part1}, there holds
\begin{align}
\sum_{i=0}^{p-1}\sum_{j=0}^{p-1}\frac{p^2(-1)^{i+j}}{(2i-3)(2j-3)(i+j+1)}
=\frac{1}{4}\left(\sum_{i=0}^{p-1}\frac{p(-1)^{i}}{2i-3}\right)^2
-\frac{1}{4}\sum_{i=0}^{p-1}\sum_{j=0}^{p-1}\frac{p^2(-1)^{i+j}}{(2i-3)(i+j+1)}. \label{eq:2i-3-first}
\end{align}
By \eqref{eq:sum-part200}, we get
\begin{align}
\sum_{i=0}^{p-1}\frac{p(-1)^{i}}{2i-3}
&=-\frac{p}{3}+\frac{p}{2p-3}-\sum_{i=0}^{p-1}\frac{p(-1)^{i}}{2i-1} \notag \\
&\equiv \frac{5p}{3}+\frac{p}{2p-3}+(-1)^{\frac{p-1}{2}}-\frac{1}{2}\sum_{k=1}^{\frac{p-3}{2}}(-1)^{\frac{p-1}{2}+k}\frac{p^2}{k^2}
\pmod{p^3},  \label{eq:2i-3}
\end{align}
It follows that
\begin{align}
\left(\sum_{i=0}^{p-1}\frac{p(-1)^{i}}{2i-3}\right)^2
\equiv 1+\frac{16p^2}{9}+\frac{10p(-1)^{\frac{p-1}{2}}}{3}+\frac{2p(-1)^{\frac{p-1}{2}}}{2p-3}
-p^2\sum_{k=1}^{\frac{p-3}{2}}\frac{(-1)^k}{k^2}  \pmod{p^3}. \label{eq:2i-3-square}
\end{align}
Moreover, we have
\begin{align*}
&\hskip -2mm \sum_{i=0}^{p-1}\sum_{j=0}^{p-1}\frac{p^2(-1)^{i+j}}{(2i-3)(i+j+1)}  \\
&=\sum_{i=0}^{p-1}\sum_{j=0}^{p-1}\frac{p^2(-1)^{i+j}}{(2i-1)(i+j+1)}
+\sum_{j=1}^{p-1}\left(\frac{p^2(-1)^j}{3j}
-\frac{p^2(-1)^j}{(2p-3)(j+p)}\right) \\
&\quad{}+\sum_{i=0}^{p-1}\left(\frac{p^2(-1)^i}{i(2i-3)}+\frac{p^2(-1)^i}{(2i-3)(i+p)}\right),
\end{align*}
and
\begin{align*}
\sum_{j=0}^{p-1}\frac{p^2(-1)^j}{(2p-3)(j+p)}
&\equiv-\sum_{j=1}^{p-1}\frac{p^2(-1)^j}{3j}+\frac{p}{2p-3} \pmod{p^3}, \\
\sum_{i=1}^{p-1}\frac{p^2(-1)^i}{i(2i-3)}
&=\frac{2}{3}\sum_{i=1}^{p-1}\frac{p^2(-1)^i}{2i-3}-\sum_{i=1}^{p-1}\frac{p^2(-1)^i}{3i}  \\
&\equiv\frac{10p^2}{9}+\frac{2p(-1)^{\frac{p-1}{2}}}{3}-\sum_{i=1}^{p-1}\frac{p^2(-1)^i}{3i} \pmod{p^3}\quad\text{(by \eqref{eq:2i-3})},\\
\sum_{i=0}^{p-1}\frac{p^2(-1)^i}{(2i-3)(i+p)}
&=\frac{2}{2p+3}\sum_{i=0}^{p-1}\frac{p^2(-1)^i}{2i-3}-\frac{1}{2p+3}\sum_{i=0}^{p-1}\frac{p^2(-1)^i}{i+p}\\
&\equiv\frac{8p^2}{9}+\frac{2p(-1)^{\frac{p-1}{2}}}{2p+3}-\frac{p}{2p+3}-\sum_{i=1}^{p-1}\frac{p^2(-1)^i}{3i}
\pmod{p^3}\quad\text{(by \eqref{eq:2i-3})}.
\end{align*}
Hence,
\begin{align}
&\hskip -2mm \sum_{i=0}^{p-1}\sum_{j=0}^{p-1}\frac{p^2(-1)^{i+j}}{(2i-3)(i+j+1)}  \notag \\
&\equiv \sum_{i=0}^{p-1}\sum_{j=0}^{p-1}\frac{p^2(-1)^{i+j}}{(2i-1)(i+j+1)}
+2p^2+\frac{2p(-1)^{\frac{p-1}{2}}}{3}+\frac{2p(-1)^{\frac{p-1}{2}}}{2p+3}-\frac{4p^2}{4p^2-9} \notag\\
&\equiv 1+\frac{4p^2}{9}+\frac{2p(-1)^{\frac{p-1}{2}}}{2p+3}
-\frac{10p(-1)^{\frac{p-1}{2}}}{3}-p^2\sum_{k=1}^{\frac{p-3}{2}}\frac{(-1)^k}{k^2}\quad\text{(by \eqref{eq:sum-part3})}.
\label{eq:2i-3-another}
\end{align}
Substituting \eqref{eq:2i-3-square} and \eqref{eq:2i-3-another} into \eqref{eq:2i-3-first}, and noticing that
\begin{align*}
\frac{p}{2p-3}-\frac{p}{2p+3}\equiv-\frac{2p}{3}\pmod{p^3},
\end{align*}
we are led to \eqref{eq:p2i-3}. \qed

\section{Proof of Theorem \ref{thm:2}}
\noindent{\it Proof of \eqref{eq:5}.}
We start with the following identity (see \cite[(2.3)]{Guo} or the proof of \cite[Lemma 4.2]{GZ2}):
\begin{align}
{k\choose i}{k+i\choose i}{k\choose j}{k+j\choose j}
&=\sum_{r=0}^i {i+j\choose i}{j\choose i-r}{j+r\choose r}{k\choose j+r}{k+j+r\choose j+r} \notag\\
&=\sum_{s=j}^{i+j} {i+j\choose i}{j\choose s-i}{s\choose j}{k\choose s}{k+s\choose s}.  \label{eq:from-GZ}
\end{align}
Noticing that ${k\choose i}{k+i\choose i}={k+i\choose 2i}{2i\choose i}$, by \eqref{eq:from-GZ}, we have
\begin{align}
&\sum_{k=0}^{n-1}(2k+1)R_k^2   \nonumber\\
&=\sum_{k=0}^{n-1}(2k+1) \sum_{i=0}^{k}\sum_{j=0}^{k}
{k+i\choose 2i}{2i\choose i}{k+j\choose 2j}{2j\choose j}\frac{1}{(2i-1)(2j-1)}  \nonumber\\
&=\sum_{k=0}^{n-1}(2k+1) \sum_{i=0}^{k}\sum_{j=0}^{k}\sum_{s=j}^{i+j}
{i+j\choose i}{j\choose s-i}{s\choose j}{k\choose s}{k+s\choose s}\frac{1}{(2i-1)(2j-1)}  \nonumber\\
&=\sum_{k=0}^{n-1}(2k+1)\sum_{m=0}^{2k}\sum_{s=0}^m
\sum_{i=0}^s {m\choose i}{m-i\choose m-s}{s\choose m-i}{k\choose s}{k+s\choose s}\frac{1}{(2i-1)(2m-2i-1)},   \label{eq:step01}
\end{align}
where $m=i+j$.

Noticing that ${m-i\choose m-s}{s\choose m-i}={s\choose i}{i\choose m-s}$ and the easily checked identity
\begin{align*}
\sum_{k=s}^{n-1}(2k+1){k\choose s}{k+s\choose s}=n{n+s\choose 2s}{2s\choose s}\frac{n-s}{s+1},
\end{align*}
we may simplify \eqref{eq:step01} as
\begin{align}
&\hskip -2mm
\sum_{k=0}^{n-1}(2k+1)R_k^2   \nonumber\\
&=n\sum_{s=0}^{n-1}{n+s\choose 2s}{2s\choose s}\frac{n-s}{s+1}\sum_{m=s}^{2s}
\sum_{i=0}^s {m\choose i}{s\choose i}{i\choose m-s}\frac{1}{(2i-1)(2m-2i-1)}.  \label{eq:step1}
\end{align}
It follows from \eqref{eq:10} that the right-hand side of \eqref{eq:step1} is an integer divisible by $n$.
This completes the proof of \eqref{eq:5}.  \qed

\medskip
\noindent{\it Proof of \eqref{eq:new5}.}
Let $n=p$ be an odd prime in \eqref{eq:step1}. For $0\leqslant s\leqslant p-2$, we have
$$
{p+s\choose 2s}{2s\choose s}\frac{p-s}{s+1}
={p-1\choose s}{p+s\choose s}\frac{p}{s+1}=\frac{p}{s+1}\prod_{i=1}^s\frac{p^2-i^2}{i^2}\equiv \frac{p(-1)^s}{s+1} \pmod{p^3}.
$$
Furthermore, by \eqref{eq:10}, for $0\leqslant i\leqslant s$ and $s\leqslant m\leqslant 2s$,
$$
{2s\choose s}{m\choose i}{s\choose i}{i\choose m-s}\frac{1}{(2i-1)(2m-2i-1)}
$$
is an integer, and since ${2s\choose s}\not\equiv 0\pmod{p^2}$, we conclude that
$$
{m\choose i}{s\choose i}{i\choose m-s}\frac{p}{(2i-1)(2m-2i-1)}
$$
is a $p$-adic integer. Similarly, for $s=p-1$, we have
$$
{p+s\choose 2s}{2s\choose s}\frac{p-s}{s+1}\equiv \frac{p(-1)^s}{s+1}=(-1)^s \pmod{p^2},
$$
and  by \eqref{eq:10}, for $0\leqslant i\leqslant p-1$ and $p-1\leqslant m\leqslant 2p-2$,
$$
{2p-2\choose p-1}{m\choose i}{p-1\choose i}{i\choose m-p+1}\frac{1}{(2i-1)(2m-2i-1)}\equiv 0\pmod{p},
$$
and so
$$
{m\choose i}{p-1\choose i}{i\choose m-p+1}\frac{p}{(2i-1)(2m-2i-1)}\equiv 0\pmod{p}.
$$

It follows that
\begin{align}
&\hskip -2mm
\sum_{k=0}^{p-1}(2k+1)R_k^2   \nonumber\\
&\equiv \sum_{s=0}^{p-1}\frac{p(-1)^{s}}{s+1}\sum_{m=s}^{2s}
\sum_{i=0}^s {m\choose i}{s\choose i}{i\choose m-s}\frac{p}{(2i-1)(2m-2i-1)}  \nonumber  \\
&=\sum_{m=0}^{p-1}\sum_{i=0}^{m}\sum_{s=i}^{m}\frac{p^2(-1)^{s}}{s+1}
{m\choose i}{s\choose i}{i\choose m-s}\frac{1}{(2i-1)(2m-2i-1)}  \nonumber \\
&\quad{}+\sum_{m=p}^{2p-2}\sum_{i=m-p+1}^{p-1}\sum_{s=i}^{p-1}\frac{p^2(-1)^{s}}{s+1}
{m\choose i}{s\choose i}{i\choose m-s}\frac{1}{(2i-1)(2m-2i-1)}.  \pmod{p^3} \label{eq:step2}
\end{align}
Note that, for $p\leqslant m\leqslant 2p-2$, $m-p+1\leqslant i\leqslant p-1$, and $p\leqslant s\leqslant m$,
we have ${m\choose i}\equiv{s\choose i}\equiv 0\pmod{p}$, $s+1\not\equiv 0\pmod{p}$. Thus, it is not hard to see that
\begin{align}
\frac{p^2(-1)^{s}}{s+1}{m\choose i}{s\choose i}{i\choose m-s}\frac{1}{(2i-1)(2m-2i-1)}
\equiv 0\pmod{p^3}, \label{eq:cases}
\end{align}
except for $m=p+1$ and $i=\frac{p+1}{2}$ ($s=p$ or $s=p+1$), in which case the left-hand side of \eqref{eq:cases}
is equal to
$$
\begin{cases}
\displaystyle\frac{-1}{2}{p+1\choose \frac{p+1}{2}}{p\choose \frac{p+1}{2}},&\text{if $s=p$,}\\[10pt]
\displaystyle\frac{1}{p+2}{p+1\choose \frac{p+1}{2}}^2,&\text{if $s=p+1$.}
\end{cases}
$$
Therefore, using the following identity
\begin{align*}
\sum_{s=i}^{m}\frac{(-1)^s}{s+1}{s\choose i}{i\choose m-s}=\frac{(-1)^m}{m+1}{m\choose i}^{-1},
\end{align*}
which can be proved by the Zeilberger algorithm (see \cite{Koepf,PWZ}), we deduce from \eqref{eq:step2} that
\begin{align*}
&\hskip -2mm \sum_{k=0}^{p-1}(2k+1)R_k^2  \\
&\equiv \sum_{m=0}^{p-1}(-1)^m\sum_{i=0}^{m} \frac{p^2}{(2i-1)(2m-2i-1)(m+1)}  \\
&\quad{}+\sum_{m=p}^{2p-2}(-1)^m\sum_{i=m-p+1}^{p-1} \frac{p^2}{(2i-1)(2m-2i-1)(m+1)}  \\
&\quad{}+\frac{1}{2}{p+1\choose \frac{p+1}{2}}{p\choose \frac{p+1}{2}}-\frac{1}{p+2}{p+1\choose \frac{p+1}{2}}^2 \\
&=\sum_{i=0}^{p-1}\sum_{j=0}^{p-1}\frac{p^2(-1)^{i+j}}{(2i-1)(2j-1)(i+j+1)}+\frac{4(p-2)p^2}{(p+2)(p-1)^2}{p-1\choose \frac{p-1}{2}}^2   \pmod{p^3}.
\end{align*}
Applying \eqref{eq:p2i-1} and noticing that
\begin{align}
\frac{4(p-2)p^2}{(p+2)(p-1)^2}{p-1\choose \frac{p-1}{2}}^2
\equiv -4p^2 \pmod{p^3},
\end{align}
we complete the proof of \eqref{eq:new5}.  \qed

\section{Proof of Theorem \ref{thm:3}}
\noindent{\it Proof of \eqref{eq:6}.} Similarly to \eqref{eq:step1}, we have
\begin{align*}
&\hskip -2mm
9\sum_{k=0}^{n-1}(2k+1)W_k^2   \nonumber\\
&=n\sum_{s=0}^{n-1}{n+s\choose 2s}{2s\choose s}\frac{n-s}{s+1}\sum_{m=s}^{2s}
\sum_{i=0}^s {m\choose i}{s\choose i}{i\choose m-s}\frac{81}{(2i-3)(2m-2i-3)}.
\end{align*}
The proof then follows from \eqref{eq:16}.  \qed

\medskip
\noindent{\it Proof of \eqref{eq:new6}.}
For $p\leqslant m\leqslant 2p-2$, $m-p+1\leqslant i\leqslant p-1$, and $p\leqslant s\leqslant m$,
we have
\begin{align}
\frac{p^2(-1)^{s}}{s+1}{m\choose i}{s\choose i}{i\choose m-s}\frac{1}{(2i-3)(2m-2i-3)}
\equiv 0\pmod{p^3}, \label{eq:cases2}
\end{align}
except for $m=p+3$ and $i=\frac{p+3}{2}$ ($s=p,\ldots,p+3$), in which case the left-hand side of \eqref{eq:cases2}
is equal to
$$
\begin{cases}
\displaystyle\frac{-(p+3)(p-1)}{48}{p+3\choose \frac{p+3}{2}}{p\choose \frac{p+3}{2}},&\text{if $s=p$,}\\[10pt]
\displaystyle\frac{(p+3)(p+1)}{8(p+2)}{p+3\choose \frac{p+3}{2}}{p+1\choose \frac{p+3}{2}},&\text{if $s=p+1$,}\\[10pt]
\displaystyle\frac{-1}{2}{p+3\choose \frac{p+3}{2}}{p+2\choose \frac{p+3}{2}},&\text{if $s=p+2$,}\\[10pt]
\displaystyle\frac{1}{p+4}{p+3\choose \frac{p+3}{2}}^2,&\text{if $s=p+3$.}
\end{cases}
$$
Similarly to the proof of \eqref{eq:new5}, we have
\begin{align*}
&\hskip -2mm \frac{1}{9}\sum_{k=0}^{p-1}(2k+1)W_k^2  \\
&=\sum_{i=0}^{p-1}\sum_{j=0}^{p-1}\frac{p^2(-1)^{i+j}}{(2i-3)(2j-3)(i+j+1)}
+\frac{(p+3)(p-1)}{48}{p+3\choose \frac{p+3}{2}}{p\choose \frac{p+3}{2}}    \\
&\quad{}-\frac{(p+3)(p+1)}{8(p+2)}{p+3\choose \frac{p+3}{2}}{p+1\choose \frac{p+3}{2}}
+\frac{1}{2}{p+3\choose \frac{p+3}{2}}{p+2\choose \frac{p+3}{2}}
-\frac{1}{p+4}{p+3\choose \frac{p+3}{2}}^2 \\
&=\sum_{i=0}^{p-1}\sum_{j=0}^{p-1}\frac{p^2(-1)^{i+j}}{(2i-3)(2j-3)(i+j+1)} \\
&\quad{}+\frac{2p^2(p^5-5p^4-3p^3+41p^2-10p-120)}{3(p+1)^2(p+3)^2(p+4))}{p-1\choose\frac{p-1}{2}}^2 \\
&\equiv \frac{4p}{3}(-1)^{\frac{p-1}{2}}+\frac{p^2}{3}-\frac{20}{9}p^2 \pmod{p^3} \quad\text{(by \eqref{eq:p2i-3})}.
\end{align*}
This completes the proof. \qed

\section{Some open problems}
Our proof of Lemmas \ref{lem:3} and \ref{lem:4} depends heavily on the multi-variable Zeilberger algorithm. Even worse,
the recurrences produced by this algorithm are very complicated. It is natural to ask the following question:
\begin{prob}{\rm
Is there any simple proof of Lemmas \ref{lem:3} and \ref{lem:4}? }
\end{prob}

It seems that the congruences \eqref{eq:5} and \eqref{eq:6} have the following refinement.
\begin{conj}
Let $n$ be a positive integer. Then
\begin{align*}
\sum_{k=0}^{n-1}(2k+1)R_k^2 &\equiv n^2 \pmod{16n}, \\
\sum_{k=0}^{n-1}(2k+1)W_k^2 &\equiv n^2 \pmod{8n}.
\end{align*}
\end{conj}

Let
$$
R_{n,r}=\sum_{k=0}^{n}{n+k\choose 2k}{2k\choose k}\frac{1}{2k-2r-1}.
$$
Then numerical calculation suggests the following conjecture.
\begin{conj}
Let $n$ and $r$ be positive integers. Then there exists an integer $a_{r}$, independent of $n$, such that
\begin{align}
a_{r}\sum_{k=0}^{n-1}(2k+1)R_{k,r}^2\equiv 0 \pmod{n}. \label{eq:22}
\end{align}
\end{conj}

We also think that $a_r=(2r+1)!!^2=(2r+1)^2(2r-1)^2\cdots 3^2\cdot 1$ is a suitable choice for \eqref{eq:22}.

\vskip 5mm \noindent{\bf Acknowledgments.} The authors would like to
thank Professor Michael Stoll for helpful comments through MathOverflow.
This work was partially supported by the Fundamental Research Funds for the Central Universities and the National Natural
Science Foundation of China (grant 11371144).


\begin{thebibliography}{99}
\small \setlength{\itemsep}{-.8mm}

\bibitem{AZ}M. Apagodu and D. Zeilberger, Multi-variable Zeilberger and Almkvist--Zeilberger
algorithms and the sharpening of Wilf--Zeilberger theory, Adv. Appl. Math. 37 (2006), 139--152.

\bibitem{Bober}J.W. Bober, Factorial ratios, hypergeometric series, and a family of step functions,
J. London Math. Soc. 79 (2009), 422-444.

\bibitem{Guo2013}V.J.W. Guo, Proof of Sun's conjecture on the divisibility of certain binomial sums,
Electron. J. Combin. 20(4) (2013), \#P20.

\bibitem{Guo2014}V.J.W. Guo, Proof of two divisibility properties of binomial coefficients conjectured by Z.-W. Sun,
Electron. J. Combin. 21(2) (2014), \#P2.54.

\bibitem{Guo}V.J.W. Guo, Some congruences involving powers of Legendre polynomials, preprint, 2014, arXiv:1412.7724.

\bibitem{GK}V.J.W. Guo and C. Krattenthaler, Some divisibility properties of binomial and $q$-binomial coefficients,
J. Number Theory 135 (2014), 167--184.

\bibitem{GZ2}V.J.W. Guo and J. Zeng, Proof of some conjectures of Z.-W. Sun on congruences for Ap\'ery polynomials,
J. Number Theory, 132 (2012), 1731--1740.

\bibitem{Koepf}W. Koepf, Hypergeometric Summation, an Algorithmic Approach to Summation
and Special Function Identities, Friedr. Vieweg \& Sohn, Braunschweig, 1998.

\bibitem{PWZ}M. Petkov\v{s}ek, H. S. Wilf and D. Zeilberger, $A=B$, A K Peters, Ltd., Wellesley, MA, 1996.

\bibitem{Sun2012}Z.-W. Sun, On divisibility of binomial coefficients, J. Austral. Math. Soc. 93 (2012), 189--201.

\bibitem{Sun2013}Z.-W. Sun, Products and sums divisible by central binomial coefficients, Electron. J. Combin. 20(1) (2013), \#P9.

\bibitem{Sun}Z.-W. Sun, Two new kinds of numbers and related divisibility results, preprint, 2014, arXiv:1408.5381v8.

\bibitem{WZ}H.S. Wilf and D. Zeilberger, An algorithmic proof theory for hypergeometric (ordinary and ``$q$") multisum/integral identities,
Invent. Math. 108 (1992), 575--633.

\end{thebibliography}
\end{document}